\documentclass{amsart}
 \newtheorem{thm}{Theorem}[section]
 \newtheorem{cor}[thm]{Corollary}
 
 \newtheorem{prop}[thm]{Proposition}
 \theoremstyle{definition}
 \newtheorem{defn}[thm]{Definition}
 \theoremstyle{remark}
 \newtheorem{rem}[thm]{Remark}
 
 \numberwithin{equation}{section}

\begin{document}
%
%
%
%
%
%
%
%
%
\title[First Order Algebraic Differential Equations]
 {Algorithmic Reduction and Rational General Solutions of First Order Algebraic
  Differential Equations}
\author[Guoting CHEN]{Guoting CHEN }

\address{Laboratoire P. Painlev\'e \\
UMR CNRS 8524 \\
UFR de Math\'{e}matiques \\
Universit\'{e} de Lille 1 \\
59655 Villeneuve d'Ascq, France }

\email{Guoting.Chen@math.univ-lille1.fr}

\thanks{The second  author is partially supported by the NKBRSF of China
(No. 2004CB318000), the NNSF (No. 10301032) and by a CNRS---K. C.
WONG fellowship during his visit to the Laboratoire P. Painlev\'e,
Universit\'e de Lille 1, France. }
\author[Yujie MA]{Yujie MA }
\address{Key Laboratory of Mathematics Mechanization\\
Chinese Academy of Sciences\\
 Beijing 100080, P. R. China } \email{yjma@mmrc.iss.ac.cn}
\subjclass{Primary 34A09, 68W30; Secondary 14Q05, 34M15}

\keywords{Algebraic differential equation, rational general
solution, rational parametrization}

\date{January 1, 2004}

\begin{abstract}
 First order algebraic differential
equations are considered. An necessary condition for a first order
algebraic differential equation to have a rational general
solution is given: the algebraic genus of the equation should be
zero. Combining with Fuchs' conditions\index{Fuchs's conditions}
for algebraic differential equations without movable critical
point, an algorithm is given for the computation of rational
general solutions of these equations if they exist under the
assumption that a rational parametrization is provided. It is
based on an algorithmic reduction of first order algebraic
differential equations with algebraic genus zero and without
movable critical point to classical Riccati
equations\index{Riccati equation}.
\end{abstract}

\maketitle

\section{Introduction}

The study of first order algebraic differential
equations\index{first order algebraic differential equation} can
be dated back to C. Briot\index{Briot} and T.
Bouquet\index{Bouquet} \cite{Brio-Bouq}, L. Fuchs\index{Fuchs}
\cite{Fuchs} and H. Poincar\'e\index{Poincar\'e} \cite{Poin1885}.
M. Matsuda\index{Matsuda} \cite{Matsuda80} gave a modern
interpretation of the results using the theory of differential
algebraic function field of one variable, and
Eremenko\index{Eremenko} \cite{eremenko98} studied the bound of
the degrees of the rational solutions of  a first order algebraic
differential equation by using the approach of \cite{Matsuda80}.

From an algorithmic point of view, many authors have been
interested in the constructions of closed form solutions for
differential equations (this problem can be traced back to the
work of Liouville\index{Liouville}). In \cite{risch},
Risch\index{Risch}, gave an algorithm for finding closed forms for
integration. In \cite{kov}, Kovacic\index{Kovacic} presented a
method for solving second order linear homogeneous differential
equations. In \cite{singer}, Singer\index{Singer} proposed a
method for finding liouvillian solutions of general linear
differential equations. In \cite{lizm}, Li\index{Li} and
Schwarz\index{Schwarz} gave a method to find rational solutions
for a class of partial differential equations. All these works are
limited to linear cases.

For algebraic (nonlinear) differential equations there are some
studies in this direction. For Riccati equations\index{Riccati
equation}, polynomial solutions are considered  in
\cite{Campbell52} and algorithms for the computation of rational
solutions are given in \cite{Bron92,kov}. In \cite{carnicer}, the
algebraic solutions of a general class of first order and first
degree algebraic differential equations was studied and the degree
bound of algebraic solutions in the nondicritical case was given.
In \cite{feng-gao, ma-gao}, algorithms for the computation of
rational general solutions or polynomial solutions are given for
some kinds of algebraic differential equations.

Another motivation of our work is differential algebraic
 geometry. In the series papers of Wu\index{Wu} \cite{wu87, wu89-1, wu03}, the author studied algebraic
differential geometry from several different points of view. In
\cite{wu03}, the author presents an algorithmic method of solving
arbitrary systems of algebraic differential equations by extending
the characteristic set method to differential case. The Devil's
problem of Pommaret\index{Pommaret} is solved in details as an
illustration.

In this paper, we consider the computation of rational general
solutions of first order algebraic differential equations by using
methods from algebraic geometry. We give an
  necessary condition for a first order algebraic differential
equation to have a  rational general solution: that is the
algebraic genus of the equation should be zero. Combining with
Fuchs' conditions\index{Fuchs' conditions} for first order
algebraic differential equations without movable critical
points\index{movable critical point}, we obtain an algorithm for
the computation of rational general solutions under the assumption
that a rational parametrization is provided. It is based on an
algorithmic reduction of first order algebraic differential
equation of algebraic genus zero and without movable critical
point to classical Riccati equations\index{Riccati equation}.

\section{Rational General Solutions of First Order Algebraic Differential
Equations }

We first present some results from algebraic geometry, which is
used in the following.

Let $f(x,y)$ be an irreducible polynomial over $\mathbb{C}$. We
call that $f(x,y)=0$ is a rational curve if there exist two
rational functions $\phi(t)$, $\psi(t)\in \mathbb{C}(t)$ such that
\begin{itemize}
\item[(i)] For all but a finite set of $t_0\in \mathbb{C}$,
$(\phi(t), \psi(t))$ is a point of $f$. \item[(ii)] With a finite
number of exceptions, for every point $(x_0, y_0)$ of $f$ there is
a unique $t_0\in \mathbb{C}$ such that $x_0=\phi(t_0)$,
$y_0=\psi(t_0)$.
\end{itemize}

It is impossible to avoid having a finite number of exceptions in
the above conditions. They arise from two sources. One is the fact
that a rational function is not defined for some values of the
variable, and the other is the presence of singular points on the
curve.

The following results are well known in algebraic geometry
\cite{walker}.

\begin{thm}\label{rational_curve}
An algebraic curve is rational if and only if its genus is zero.
\end{thm}

\begin{thm}\label{luroth}
(a) Every rational transform of a rational curve is a rational
curve.

\noindent (b) If $\lambda$ is transcendental over $\mathbb{C}$ and
if $\mathbb{C} \subset F \subset \mathbb{C}(\lambda)$, $F \not=
\mathbb{C}$, then there is an element $\mu\in F$, transcendental
over $\mathbb{C}$, such that $F = \mathbb{C}(\mu)$.

\noindent (c) If a curve $f(x,y)=0$ satisfies (i) for rational
functions $\phi(\lambda), \psi(\lambda)$ which are not both
constants, then there exist rational functions
$\tilde\phi(\lambda), \tilde\psi(\lambda)$ for which both (i) and
(ii) are satisfied, and the curve is rational.
\end{thm}

The three statements in Theorem \ref{luroth} are all equivalent.
It is called L\"uroth's theorem\index{L\"uroth's theorem}.

Now we consider a first order algebraic differential equation in
the form
\begin{equation}\label{alg_diff_eq}
F(z,w,w')=0,
\end{equation}
where $F$ is a polynomial in $w,w'$ with rational coefficients in
$\mathbb{C}(z)$.

\begin{defn}A general solution $w(z,\lambda)$ of (\ref{alg_diff_eq}) is called a rational general solution\index{rational general solution} if it
is rational in $z$ and $\lambda$.
\end{defn}

If $F(z,w,w^{\prime})=0$ admits a rational general solution, then
it is free from movable critical point for poles are the only
singularities of the solution which change their position if one
varies the initial data $c\in \mathbb{C}$.

We now prove the following theorem on rational general solution of
a first order algebraic differential equation.

\begin{thm}\label{main_th}
If a first order irreducible algebraic differential equation
$$F(z,w,w^{\prime})=0$$
admits a non-constant rational general solution, then the genus of
$F(z,w,w^{\prime})=0$ with respect to $w,w'$ is zero for any $z$,
except for a finitely many exceptions.
\end{thm}
\begin{proof} Let $w=r(z,\lambda)$ be the rational general solution of
$F=0$ with the arbitrary constant $\lambda$.  Then
$w(z)=r(z,\lambda)$ and $w'(z)=\frac{\partial r}{\partial z}$ are
rational functions and they satisfy equation $F(z,w,w')=0$.

Let $z$ be fixed and consider the curve $f_z(x,y)=F(z,x,y)=0$.
Denote $\phi_z(\lambda)= w(z)$ and $\psi_z(\lambda)= w'(z)$. If
$\psi_z(\lambda)$ is a constant, then
$w(z)=zw^{\prime}(z)+\lambda$ is of genus zero and hence the genus
of $f_z(x,y)=0$ is also zero. If $\psi_z(\lambda)$ is not a
constant, consider the point $(\phi_z(\lambda),\psi_z(\lambda))$
of $f_z(x,y)=0$ for the parameter $\lambda$ in the transcendental
extension field $\mathbb{C}(\lambda)$. It is clear that (i) is
satisfied for all but finitely many $\lambda$. Hence $f_z(x,y)=0$
is a rational curve and its genus is zero by Theorem
\ref{rational_curve} and Theorem \ref{luroth}.
\end{proof}

Motivated by this theorem, we present the following definition.

\begin{defn}The algebraic genus\index{algebraic genus} of a first order algebraic
differential equation $F(z, w, w^{\prime})=0$ is defined to be the
genus of $F(z, w, w^{\prime})=0$ with respect to $w$ and
$w^{\prime}$.
\end{defn}

\section{Reduction of First Order Algebraic Differential Equations} \label{sec_redu}

For a first order algebraic differential equation
$F(z,w,w^{\prime})=0$, Fuchs Theorem\index{Fuchs Theorem}
presented necessary conditions for the equation to be free from
movable critical point.  By the Painlev\'{e}
Theorem\index{Painlev\'{e} Theorem}, we know that Fuchs'
conditions\index{Fuchs' conditions} are sufficient (see
\cite{Fuchs, ince, Pain,Poin1885}).

\subsection{Fuchs Theorem\index{Fuchs Theorem}}

Let $D(z, w)$ be the $p$--discriminant of the equation
$F(z,w,w^{\prime})=0,$ it is a polynomial in $w$, whose
coefficients are analytic functions of $z$ \cite{ince}.

The conditions, necessary  to secure that the first order
differential equation
$$F(z,w,w^{\prime})=0$$ of degree $m$, shall have no movable critical point,
are:
\begin{enumerate}
\item The coefficient $A_0(z, w)$ is independent of $w$ and
therefore reduces to a function of $z$ alone or to a constant. The
equation may then be divided throughout by $A_0$ and takes the
form $$w^{\prime m}+\psi_1(z, w)w^{\prime m-1}+\cdots
+\psi_{m-1}(z, w)w^{\prime}+\psi_m(z, w)=0$$ in which the
coefficients $\psi$ are polynomials in $w$, and analytic, except
for isolated singular points, in $z$.

\item If $w=\eta(z)$ is a root of $D(z, w)=0$, and $p=\omega(z)$
is a multiple root of $F(z, \eta, \eta^{\prime})=0$, such that the
corresponding root of $F(z, w, w^{\prime})=0$, regarded as a
function of $w-\eta(z)$ is branched, then
$$\omega(z)=\frac{d\eta}{dz}.$$

\item If the order of any branch is $\alpha$, so that the equation
is effectively of the form
$$\frac{d}{dz}\{w-\eta(z)\}=c_k\{w-\eta(z)\}^{\frac{k}{\alpha}}$$ then $k\geq\alpha-1$.
\end{enumerate}

\subsection{Reduction to Classical Riccati Equation\index{Riccati equation}}\label{red_riccati}

Consider now a first order algebraic differential equation
$F(z,w,w')=0$ of genus zero and without movable critical point.
One can find a parametrization of the rational curve $F(z,x,y)=0$
in the form $x=r_1(t,z)$ and $y=r_2(t,z)$ with $r_1(t,z)$ and
$r_2(t, z)$ rational functions in $t$ and $z$. By the inversion of
rational curves we know that $t$ is rational function in $z$, $x$
and $y$. For algorithm on parametrization and inversion of
rational curves we refer to \cite{abh-baj,
sederberg-anderson-goldman, sederberg-zheng, Sen-Win, Sen-Win2,
vanH}. One has
\begin{equation}\label{ricc_eq1}
\frac{dt}{dz}=\Big(r_2(t,z)-\frac{\partial r_1}{\partial
z}\Big)/\frac{\partial r_1}{\partial t} = \frac{P(t,z)}{Q(t,z)},
\end{equation}
where $P$ and $Q$ are polynomials in $t$ and $z$. Since
$F(z,w,w')=0$ has no movable critical point, one knows that
$(\ref{ricc_eq1})$ is also has no movable critical point as $t$ is
rational function in $z$, $x$ and $y$. By Fuchs
Theorem\index{Fuchs Theorem}, we have that equation
$(\ref{ricc_eq1})$ is a Riccati equation\index{Riccati equation}
(\cite{golubev1} Chapter II, \S 7), that is
\begin{equation}\label{ricc_eq2}
\frac{dt}{dz}= A(z)t^2+B(z)t+C(z),
\end{equation}
where $A,B,C$ are rational functions in $z$. We distinguish two
cases according to $A(z)$.

\vspace{3mm}

\noindent {\bf Case 1:} If $A(z)\not\equiv0$, we consider the
change of variables $t(z)=-u(z)/A(z)$. One has
$$u'(z) +u^2 = (B(z) + A'(z)/A(z)) u -C(z)A(z),$$
In which the coefficient $A(z)$ is reduced to $-1$. Next we make
the change $u=v+\beta(z)$ to reduce the coefficient of $u$ to zero
by choosing an appropriate $\beta$. We obtain finally a classical
Riccati equation in the form
\begin{equation} \label{classical_Ricc}
v'+v^2 = r(z) \in \mathbb{C}(z).
\end{equation}
Algorithms for the computations of rational solutions of classical
Riccati equations are given in many literatures (see for example
\cite{Bron92, kov}).

If $r(z)\not\equiv0$ then a rational solution of equation
$(\ref{classical_Ricc})$ is equivalent to an exponential solution
$e^{\int v(z) dz}$ of the linear differential equation
\begin{equation}\label{2nd-order}
y'' = r(z) y.
\end{equation}

\begin{prop}
If the Riccati equation\index{Riccati equation}
$(\ref{classical_Ricc})$ with $r(z)\not\equiv 0$ has a general
rational solution, then $r(z)$ has the form
$$r(z)=\sum_{i=1}^m\left(\frac{\beta_i}{(z-z_i)^2} +
\frac{\gamma_i}{(z-z_i)}\right), $$ in which $4\beta_i = n_i^2-1$
where $n_i$ is an integer $\ge 2$.
\end{prop}
\begin{proof} Suppose that $v(z)$ is a rational solution of equation
$(\ref{classical_Ricc})$. Let $z_1, \cdots, z_m$ be the poles of
$r$. According to Kovacic's algorithm\index{Kovacic's algorithm}
(see \cite{kov}), $v(z)$ should be in the form
$$v(z)=\sum_{i=1}^m\sum_{j=1}^{\nu_i} \frac{a_{ij}}{(z-z_i)^j}
+\sum_{k=1}^d\frac{1}{z-c_k}+f(z),$$ where the $\nu_i$ are known,
$a_{ij}$ are known up to a two choices each, $d$ is known, and
$f\in \mathbb{C}[z]$ is known up to two choices. Hence there may
be arbitrary parameter only in the determination of the $c_k$.

Let
$$P(z)=\prod_{k=1}^{d} (z-c_k),$$ and
$$\omega(z)=\sum_{i=1}^m\sum_{j=1}^{\nu_i} \frac{a_{ij}}{(z-z_i)^j}+f(z).$$
Then $v=P'/P+\omega$ and $y=e^{\int v}=Pe^{\int\omega}$ is a
solution of the linear differential equation (\ref{2nd-order}).
Hence $P$ is a polynomial solutions of degree $d$ of the following
linear equation:
\begin{equation}\label{lin_eq}
P''+2\omega P' + (\omega'+\omega^2-r)P =0.
\end{equation}
One can determine whether it has a general polynomial solution or
not.

Furthermore if $(\ref{classical_Ricc})$ admits a rational general
solution, then writing $r(z)=p(z)/q(z)$, according to \cite{Yuan},
one has
\begin{itemize}
\item[(a)] $\deg (p) - \deg(q) \le -2$; \item[(b)] $r(z)$ has only
double poles, hence
$$r(z)=\sum_{i=1}^m\left(\frac{\beta_i}{(z-z_i)^2} +
\frac{\gamma_i}{(z-z_i)}\right), $$ in which $4\beta_i = n_i^2-1$
where $n_i$ is an integer $\ge 2$.
\end{itemize}
\end{proof}

Therefore a possible rational solution of equation
$(\ref{classical_Ricc})$ with $r(z)$ as in the proposition should
be
$$v(z)=\sum_{i=1}^m \frac{a_{i}}{z-z_i}
+\sum_{k=1}^d\frac{1}{z-c_k},$$ where $a_i^2-a_i-\beta_i=0$. Hence
to determine a rational general solution one needs to compute
polynomial solutions of equation $(\ref{lin_eq})$ in order to
determine the $c_k$.

\vspace{3mm}

\noindent {\bf Case 2:} If $A(z)\equiv0$, then one can integrate
easily the linear equation
$$t'=B(z)t+C(z)$$
to get the general solution
$$t(z)= \Big(\int C(z)dz + \lambda \Big)e^{\int B(z) dz},$$
where $\lambda$ is an arbitrary constant. Effective algorithm is
given in \cite{risch} for integration in closed forms. One may
find rational solutions in this way. It is clear that in this case
one may get rational general solutions only if both
$\displaystyle\int C(z)dz$ and $e^{\int B(z) dz}$ are rational
functions.

\subsection{First Order Algebraic Differential Equations with Constant
  Coefficients}

As an application, we consider a first order algebraic
differential equation  with constant coefficients,
\begin{equation}\label{const_coef}
F(w,w')=0.
\end{equation}
This kind of equations was systematically studied in
\cite{Brio-Bouq} and an algorithm is given in \cite{feng-gao} for
the determination of the rational general solution of the equation
if it exists.

When using the above reduction for equation $(\ref{const_coef})$,
it is clear that $z$ is not involved in the equation, we
henceforth get a Riccati equation\index{Riccati equation} with
constant coefficients.

As above there are two cases to be considered. In case 1, if
equation $(\ref{const_coef})$ has a non constant rational solution
$w(z)$ then $w(z+\lambda)$ is a general rational solution. Since
the equation $u'+u^2=c$ for a constant $c\not=0$ does not have a
rational solution, then one can reduce equation
$(\ref{const_coef})$ to the equation  $u^{\prime}+u^2=0$ in which
case we have the general solution $$u=\frac{1}{z+\lambda}.$$

In case 2, the equation can be converted to $u'=bu+c$ with
constant $b, c$. Then $u=\lambda e^{bz} - \frac{c}{b}$ if
$b\not=0$, and  $u=cz+\lambda$ if $b=0$, where $\lambda$ is an
arbitrary constant.

Summarizing the above, we then have the following
\begin{cor}
Let $F(w, w^{\prime})=0$ be a first order irreducible algebraic
differential equation with constant coefficients. Then it has a
non-constant rational general solution if and only if it can be
reduced either to a linear equation $u^{\prime}=c$ for some
constant $c$ or to a Riccati equation\index{Riccati equation} of
the form $u^{\prime}+u^2=0$.
\end{cor}

\noindent {\bf Example 1.} Consider
$$
F(y,y')={y'}^{4} - 8\,{y'}^{3} + (6 + 24\,{y})\,{y'}^{2} + 257 +
528\,{y}^{2} - 256\,{y}^{3} - 552\,{y}.
$$
This example comes from \cite{feng-gao}. One finds by computation
in Maple\index{Maple} that its algebraic genus is zero and it has
the following rational parametrization\index{rational
parametrization}.
$$y = \frac{17}{16}-27t+\frac{2187}{2}t^2+531441 t^4, \quad
y' = 78732t^3+81t-1. $$ And the corresponding Riccati equation is
$t'(z)=\frac{1}{27}$. Hence $t=\frac{1}{27}z+\lambda$ and the
general solution of the differential equation $F(y,y')=0$ is
$$y(z) = \frac{17}{16}-27(\frac{1}{27}z+\lambda)
+\frac{2187}{2}(\frac{1}{27}z+\lambda)^2+531441
(\frac{1}{27}z+\lambda)^4,$$ where $\lambda$ is an arbitrary
constant.

\vspace{3mm}
\section{Algorithm and Example}

We can now give the following algorithm on seeking for a rational
general solution of a first order algebraic differential equation.

\vspace{3mm}

\subsection{Algorithm}

{\noindent \bf Algorithm:}

{\bf Input:} $F(z,w,w^{\prime})=0$ is a first order algebraic
differential equation.

{\bf Output:} A rational general solution\index{rational general
solution} of $F(z,w,w^{\prime})=0$ if it exists.

\begin{enumerate}
\item Determine the irreducibility of the equation.

If $F(z,w,w^{\prime})=0$ is reducible, then factorize it and go to
step $(2)$ for each branch curve of $F(z,w,w^{\prime})=0$, else go
to step $2$ directly.

\item Compute the algebraic genus\index{algebraic genus} $g$ of
$F(z,w,w^{\prime})=0$. If $g\not=0$, then the equation doesn't
admit rational general solution\index{rational general solution}
by Theorem \ref{main_th}, else go to step $3$.

\item Determine the Fuchs' conditions\index{Fuchs' conditions} of
$F(z,w,w^{\prime})=0$. If the conditions are not satisfied, the
algorithm terminates, else go to step $4$.

\item If a rational parametrization\index{rational
parametrization} of  $F(z,w,w^{\prime})=0$ in the form
$w=r_1(t,z)$ and $w^{\prime}=r_2(t,z)$ with $r_1(t,z)$ and
$r_2(t,z)$ rational functions in $t$ and $z$ is provided, go to
step $5$.

\item Compute the derivative
$\frac{dt}{dz}=(r_2(t,z)-\frac{\partial r_1}{\partial
z})/\frac{\partial r_1}{\partial t}$ which is a Riccati
equation\index{Riccati equation} of the form $(\ref{ricc_eq2})$ by
the Fuchs Theorem\index{Fuchs Theorem}.

\item Reduce the above Riccati equation to a classical Riccati
equation\index{Riccati equation} $(\ref{classical_Ricc})$ and
compute a rational solution using the algorithm in
\cite{Bron92,kov}.

\end{enumerate}

\begin{rem}
Assume that  $F(z, w, w^{\prime})=0$ has algebraic genus $0$ for
all $z$ (or for all $z$ with some finite exceptions), in general,
the rational parametrization\index{rational parametrization} that
we obtain using the the algorithm in \cite{abh-baj,Sen-Win,
Sen-Win2, vanH} are not rational functions in $z$ (it could be
appear algebraic elements over $\mathbb{C}(z)$ as $\sqrt{z}$ (see
\cite{Sen-Win2})). Our algorithm works only with for those
equation $F(z, w, w^{\prime})=0$ for which  a rational
parametrization is provided. For instance, for first order
algebraic differential equations with constant coefficients.
\end{rem}

\begin{rem}Our algorithm for the case of constant coefficients
is the same as the classical one, which is for instance described
in \cite{Picard}, Ch. IV, sec. I (page 62--65).
\end{rem}

 \subsection{Example}

\noindent {\bf Example 2.} Consider the following equation:
$$F(z,w,w')={w'}^{2}+{\frac {2w}{z}}w'-4\,zw^{3}
+{\frac {\left( 1 +12\,{z}^{2} \right) {w}^{2}}{{z}^{2}}}
-12\,{\frac {w}{z}}+ \frac{4}{z^2}.
$$
Its algebraic genus is zero. One gets the following rational
parame\-trization by Maple\index{Maple}:
$$\begin{array}{l}
\displaystyle w= r_1=\frac{t^2z^2+4t^2-6tz+1+4z^2} {4z(-z+t)^2},
\cr \displaystyle w'=r_2= -\frac{-4z^3+13tz^2+t+2t^2z^3-10t^2
z+t^3z^4+t^3z^2+4t^3} {4z^2(-z+t)^{3}}.
\end{array}$$
One obtains the following Riccati equation
$$t'=\frac {\left({z}^{2}+2\right)}{2({z}^{2}+1)}{t}^{2}
+{\frac {z}{{z}^{2}+1}}t +\frac{3}{2({z}^{2}+1)}.
$$
Continue the reduction procedure of Section \ref{red_riccati} one
obtains
$$u'+u^2=Bu+C,$$
where
$$B = -\frac{z^3}{(z^2+1)(z^2+2)} \quad \mbox{and }
\quad C=\frac{3(z^2+2)}{4(z^2+1)^2}.
$$
And finally one obtains the following classical Riccati equation
$$v'+v^2=-6(z^{2}+2)^{-2},$$
which has a rational general solution as follows
$$v(z)=
-\frac{z}{z^2+2}+\frac{1}{z-\lambda}+\frac{1}{z+2/\lambda}$$ with
$\lambda$ the arbitrary constant. One finally has by substitutions
the following solution of equation $F(z,w,w')=0$:
$$w(z)={\frac {{z}^{2}{\lambda}^{2}-2\,z{\lambda}^{3}+4\,z\lambda+4+{\lambda}
^{4}-3\,{\lambda}^{2}}{ \left( z\lambda+2-{\lambda}^{2} \right)
^{2}z} }.
$$

\section{Conclusion}
In this paper, we present an algebraic geometry approach to the
study of first order algebraic differential equations. The
necessary and sufficient conditions of the rational general
solutions of the Riccati equations\index{Riccati equation} are
studied in \cite{chen-ma-1}. In \cite{chen-ma-2}, we prove that
the degree of a rational general solution of a first order
algebraic differential equation of degree $d$ is less than or
equal to $d$. The algebraic geometry approach is also used to
obtain bound of the number of the rational solutions of a first
order algebraic differential equation with algebraic genus greater
than one. In \cite{chen-ma-3}, the algebraic general solutions of
first order algebraic differential equations were studied by using
of the birational transformations of algebraic curves, and an
algorithm was presented to get an algebraic general solution of
first order algebraic differential equations without movable
critical point if the algebraic general solution exists. It is
interesting to present an effective algorithm to find the rational
general solutions or the rational solutions of a   first order
algebraic differential equation.


\subsection*{Acknowledgment}The authors are grateful to the anonymous referee
for the valuable comments and suggestions, which improved the
manuscript, and for pointing out some inaccuracies.
\end{document}